\documentclass{amsart}
\usepackage{amsfonts,hyperref}

\newtheorem{theorem}{Theorem}
\theoremstyle{plain}

\newtheorem{corollary}{Corollary}

\newtheorem{definition}{Definition}

\newtheorem{lemma}{Lemma}

\numberwithin{equation}{section}
\theoremstyle{theorem}

\begin{document}
\title[Foliations on double-twisted products]{Foliations on double-twisted products}
\author{Andr\'e Gomes}
\address{%
\parbox{.9\textwidth}{\centering
Department of Mathematics\\
Institute of Mathematics and Statistics\\
University of S\~ao Paulo\\
S\~ao Paulo - 05508-090 - SP\\
Brazil\\[5pt]}}
\email{gomes@ime.usp.br}

\begin{abstract}
In this paper we present a certain class of geodesic vector fields of the double-twisted product $\mathbb{R}\times_{(f,g)}\mathbb{R}$. Some examples of totally geodesic foliations are given. 
\end{abstract}
\maketitle
\vspace{-0.5cm}
\section{Introduction}
The concept of double-twisted product was introduced by Ponge and Reckziegel (see [5]) in the context of the semi-Riemannian Geometry. 

\definition Let $(M_{1}, g_{1})$ and $(M_{2}, g_{2})$ be (pseudo)Riemannian manifolds. Let
$\lambda_{i}: M_{1} \times M_{2}\longrightarrow \mathbb{R}$ ($i=1,2$) be a positive and differentiable function. Consider the canonical projections $\pi_{i}: M_{1}\times M_{2}\longrightarrow M_{i}$ for $i=1,2$. Then \textit{the double-twisted product
$M_{1}\times_{(\lambda_{1},\lambda_{2})}M_{2}$ of $(M_{1},g_{1})$ and $(M_{2},g_{2})$} is the differentiable manifold
$M_{1}\times M_{2}$ equipped with the (pseudo)Riemannian metric $g$ defined by
\begin{eqnarray*}
	g(X, Y)=\lambda_{_{1}}g_{_{1}}(d\pi_{_{1}}(X), d\pi_{_{1}}(Y))+\lambda_{_{2}}g_{_{2}}(d\pi_{_{2}}(X), d\pi_{_{2}}(Y))
\end{eqnarray*}
for all vector fields $X$ and $Y$ of $M_{1}\times M_{2}$.
\\
\\
\normalfont
This definition generalizes the R. L. Bishop's notion of an umbilic
product (denoted by $M_{1}\times_{\lambda}M_{2}$), which B. Y. Chen calls a twisted product and which is
the \textit{double} twisted product $M_{1}\times_{(1,\lambda)}M_{2}$ (see [1] and [2]). If in this
situation $\lambda$ only depends on the points of $M_{1}$, then $M_{1}\times_{\lambda}M_{2}$ is a warped
product by definition.
\\
\\
If $M$ and $N$ are complete Riemannian manifolds, then the double-twisted product $M\times_{(\lambda_{1},\lambda_{2})}N$ is not, in general, a complete Riemannian manifold. But, for our purposes, the following Lemma is enough (the proof is similar of that of Lemma 40 in [4])

\lemma Suppose that $M$ and $N$ are complete Riemannian manifolds. Given two differentiable functions $\lambda_{1},\lambda_{2}:M\times N\longrightarrow [1,\infty)$, then the double-twisted product $M\times_{(\lambda_{1},\lambda_{2})}N$ is a complete Riemannian manifold.
\proof Denote by $\pi:M\times N\longrightarrow M$ and $\sigma:M\times N\longrightarrow N$ the natural projections. Denote by $\mathcal{L}$ and $d$ the arc lenght and the distance in $M\times_{(\lambda_{1},\lambda_{2})}N$, respectively. We use the metric completeness criterion from the Hopf-Rinow theorem. Note first that if $v$ is tangent to $M\times_{(\lambda_{1},\lambda_{2})}N$ then, since $\lambda_{1}\geq1$ and $\lambda_{2}>0$, we have $\left\langle v,v\right\rangle\geq \left\langle d\pi(v),d\pi(v)\right\rangle$. Hence $\mathcal{L}(\alpha)\geq\mathcal{L}(\pi\circ\alpha)$ for any curve segment $\alpha$ in $M\times_{(\lambda_{1},\lambda_{2})}N$. Analogously, $\mathcal{L}(\alpha)\geq\mathcal{L}(\sigma\circ\alpha)$. Then we have the inequalities
\begin{eqnarray*}
\ \ d(x,y)\geq d(\pi(x),\pi(y)) \ \ \ and \ \ \ d(x,y)\geq d(\sigma(x),\sigma(y)),\ \ \ \ \ \forall x,y\in M\times_{(\lambda_{1},\lambda_{2})}N
\end{eqnarray*}
Both inequalities implies that if $\left((p_{k},q_{k})\right)_{k\in\mathbb{N}}$ is a Cauchy sequence in the double-twisted product $M\times_{(\lambda_{1},\lambda_{2})}N$, then $(p_{k})_{k\in\mathbb{N}}$ and $(q_{k})_{k\in\mathbb{N}}$ are Cauchy sequences, respectively, in $(M,g_{1})$ and $(N,g_{2})$. As $M$ and $N$ are complete, $(p_{k})_{k\in\mathbb{N}}$ and $(q_{k})_{k\in\mathbb{N}}$ are convergent and then $\left((p_{k},q_{k})\right)_{k\in\mathbb{N}}$ is convergent in $M\times_{(\lambda_{1},\lambda_{2})}N$ $\square$

\section{Preliminaries}
Some computations are necessary in order to find the examples. We denote by $\partial_{x}$ and $\partial_{y}$ the canonical vector fields of $\mathbb{R}^{2}$. The new metric in $\mathbb{R}^{2}$ is given by
\[
\left\{
      \begin{array}{ll}
      \left\langle \partial_{x},\partial_{x}\right\rangle_{(x,y)}=[f(x,y)]^{2}\\
      \\
      \left\langle \partial_{x},\partial_{y}\right\rangle_{(x,y)}=0\\
      \\
      \left\langle \partial_{y},\partial_{y}\right\rangle_{(x,y)}=[g(x,y)]^{2}
      \end{array}
\right.
\]
where $(x,y)\in\mathbb{R}^{2}$ and $f,g:\mathbb{R}^{2}\:-\!\!\!\!\longrightarrow(0,+\infty)$ are functions of class $\mathcal
{C}^{\infty}$. The real plane equipped with this (Riemannian) metric is the \textit{double twisted product of $\mathbb{R}$ and $\mathbb{R}$} and denoted by $\mathbb{R}\times_{(f,g)}\mathbb{R}$. We find now its connection $\nabla$ and curvature $\mathrm{R}$.
\\
\\
By metric compatibility we know that
\begin{eqnarray*}
	2ff_{x}\!\!\!&=&\!\!\!\partial_{x}\left\langle \partial_{x},\partial_{x}\right\rangle=2\left\langle \partial_{x},\nabla_{\partial_{x}}\partial_{x}\right\rangle
	\\
	2ff_{y}\!\!\!&=&\!\!\!\partial_{y}\left\langle \partial_{x},\partial_{x}\right\rangle=2\left\langle \partial_{x},\nabla_{\partial_{y}}\partial_{x}\right\rangle
	\\
	2gg_{x}\!\!\!&=&\!\!\!\partial_{x}\left\langle \partial_{y},\partial_{y}\right\rangle=2\left\langle \partial_{y},\nabla_{\partial_{x}}\partial_{y}\right\rangle
	\\
	2gg_{y}\!\!\!&=&\!\!\!\partial_{y}\left\langle \partial_{y},\partial_{y}\right\rangle=2\left\langle \partial_{y},\nabla_{\partial_{y}}\partial_{y}\right\rangle
	\\
  0\!\!\!&=&\!\!\!\partial_{x}\left\langle \partial_{x},\partial_{y}\right\rangle=\left\langle \nabla_{\partial_{x}}\partial_{x},\partial_{y}\right\rangle+\left\langle \partial_{x},\nabla_{\partial_{x}}\partial_{y}\right\rangle
  \\
  0\!\!\!&=&\!\!\!\partial_{y}\left\langle \partial_{x},\partial_{y}\right\rangle=\left\langle \nabla_{\partial_{y}}\partial_{x},\partial_{y}\right\rangle+\left\langle \partial_{x},\nabla_{\partial_{y}}\partial_{y}\right\rangle
  \\
  \nabla_{\partial_{x}}\partial_{y}\!\!\!&=&\!\!\!\nabla_{\partial_{y}}\partial_{x}
\end{eqnarray*}
where in the last equation we used $[\partial_{x},\partial_{y}]=0$.
\\
\\
Using the equations above, we find the Riemannian connection $\nabla$ of $\mathbb{R}\times_{(f,g)}\mathbb{R}$:
\begin{eqnarray*}
	\nabla_{\partial_{x}}\partial_{x}\!\!\!&=&\!\!\!\frac{1}{f^2}\left\langle \nabla_{\partial_{x}}\partial_{x},\partial_{x}\right\rangle\partial_{x}+\frac{1}{g^2}\left\langle \nabla_{\partial_{x}}\partial_{x},\partial_{y}\right\rangle\partial_{y}
	\\
	&=&\!\!\!\frac{f_{x}}{f}\partial_{x}-\frac{ff_{y}}{g^{2}}\partial_{y}
	\\
\nabla_{\partial_{y}}\partial_{y}\!\!\!&=&\!\!\!\frac{1}{f^2}\left\langle \nabla_{\partial_{y}}\partial_{y},\partial_{x}\right\rangle\partial_{x}+\frac{1}{g^2}\left\langle \nabla_{\partial_{y}}\partial_{y},\partial_{y}\right\rangle\partial_{y}
	\\
&=&\!\!\!-\frac{g_{x}g}{f^{2}}\partial_{x}+\frac{g_{y}}{g}\partial_{y}
	\\	
\nabla_{\partial_{y}}\partial_{x}\!\!\!&=&\!\!\!\frac{1}{f^2}\left\langle \nabla_{\partial_{y}}\partial_{x},\partial_{x}\right\rangle\partial_{x}+\frac{1}{g^2}\left\langle \nabla_{\partial_{y}}\partial_{x},\partial_{y}\right\rangle\partial_{y}
\\
&=&\!\!\!\frac{f_{y}}{f}\partial_{x}+\frac{g_{x}}{g}\partial_{y}
\end{eqnarray*}
Moreover, by definition we have 
\begin{eqnarray*}	
\mathrm{R}(\partial_{x},\partial_{y})\partial_{x}\!\!\!\!&=&\!\!\!\!\nabla_{\partial_{y}}\nabla_{\partial_{x}}\partial_{x}-\nabla_{\partial_{x}}\nabla_{\partial_{y}}\partial_{x}-\nabla_{[\partial_{x},\partial_{y}]}\partial_{x}
\\
&=&\!\!\!\!\nabla_{\partial_{y}}\nabla_{\partial_{x}}\partial_{x}-\nabla_{\partial_{x}}\nabla_{\partial_{y}}\partial_{x}
\\
&=&\!\!\!\!\nabla_{\partial_{y}}(\frac{f_{x}}{f}\partial_{x}-\frac{ff_{y}}{g^{2}}\partial_{y})-\nabla_{\partial_{x}}(\frac{f_{y}}{f}\partial_{x}+\frac{g_{x}}{g}\partial_{y})
\\
&=&\!\!\!\frac{f_{x}}{f}\nabla_{\partial_{y}}\partial_{x}+(\frac{f_{x}}{f})_{_{y}}\partial_{x}-\frac{ff_{y}}{g^2}\nabla_{\partial_{y}}\partial_{y}-(\frac{ff_{y}}{g^2})_{_{y}}\partial_{y}
\\ 
&&\!\!\!\!-(\frac{f_{y}}{f}\nabla_{\partial_{x}}\partial_{x}+(\frac{f_{y}}{f})_{_{x}}\partial_{x}+\frac{g_{x}}{g}\nabla_{\partial_{x}}\partial_{y}+(\frac{g_{x}}{g})_{_{x}}\partial_{y})
\\
&=&\!\!\!\!\frac{f_{x}}{f}(\frac{f_{y}}{f}\partial_{x}+\frac{g_{x}}{g}\partial_{y})+(\frac{f_{x}}{f})_{_{y}}\partial_{x}-\frac{ff_{y}}{g^2}(-\frac{g_{x}g}{f^{2}}\partial_{x}+\frac{g_{y}}{g}\partial_{y})-(\frac{ff_{y}}{g^2})_{_{y}}\partial_{y}
\\ 
&&\!\!\!\!-\left[\frac{f_{y}}{f}(\frac{f_{x}}{f}\partial_{x}-\frac{ff_{y}}{g^{2}}\partial_{y})+(\frac{f_{y}}{f})_{_{x}}\partial_{x}+\frac{g_{x}}{g}(\frac{f_{y}}{f}\partial_{x}+\frac{g_{x}}{g}\partial_{y})+(\frac{g_{x}}{g})_{_{x}}\partial_{y}\right]
\\
&=&\!\!\!\!\left[\frac{f_{y}^{2}}{g^{2}}-\frac{g_{x}^{2}}{g^{2}}+\frac{f_{x}g_{x}}{fg}-\frac{ff_{y}g_{y}}{g^{3}}-(\frac{ff_{y}}{g^{2}})_{_{y}}-(\frac{g_{x}}{g})_{_{x}}\right]\partial_{y}
\end{eqnarray*}
and then 
\begin{eqnarray*}
	\left\langle \mathrm{R}(\partial_{x},\partial_{y})\partial_{x},\partial_{y}\right\rangle=g^{2}\left[\frac{f_{y}^{2}}{g^{2}}-\frac{g_{x}^{2}}{g^{2}}+\frac{f_{x}g_{x}}{fg}-\frac{ff_{y}g_{y}}{g^{3}}-(\frac{ff_{y}}{g^{2}})_{_{y}}-(\frac{g_{x}}{g})_{_{x}}\right]
\end{eqnarray*}
Analogously, we have
\begin{eqnarray*}
\mathrm{R}(\partial_{x},\partial_{y})\partial_{y}\!\!\!\!&=&\!\!\!\!\nabla_{\partial_{y}}\nabla_{\partial_{x}}\partial_{y}-\nabla_{\partial_{x}}\nabla_{\partial_{y}}\partial_{y}-\nabla_{[\partial_{x},\partial_{y}]}\partial_{y}
\\
&=&\!\!\!\!\nabla_{\partial_{y}}\nabla_{\partial_{x}}\partial_{y}-\nabla_{\partial_{x}}\nabla_{\partial_{y}}\partial_{y}
\\
&=&\!\!\!\!\nabla_{\partial_{y}}(\frac{f_{y}}{f}\partial_{x}+\frac{g_{x}}{g}\partial_{y})-\nabla_{\partial_{x}}(-\frac{g_{x}g}{f^{2}}\partial_{x}+\frac{g_{y}}{g}\partial_{y})
\\
&=&\!\!\!\!\frac{f_{y}}{f}\nabla_{\partial_{y}}\partial_{x}+(\frac{f_{y}}{f})_{_{y}}\partial_{x}+\frac{g_{x}}{g}\nabla_{\partial_{y}}\partial_{y}+(\frac{g_{x}}{g})_{_{y}}\partial_{y}
\\ 
&&\!\!\!\!+\frac{gg_{x}}{f^{2}}\nabla_{\partial_{x}}\partial_{x}+(\frac{gg_{x}}{f^{2}})_{_{x}}\partial_{x}-\frac{g_{y}}{g}\nabla_{\partial_{x}}\partial_{y}-(\frac{g_{y}}{g})_{_{x}}\partial_{y}
\\
&=&\!\!\!\!\frac{f_{y}}{f}(\frac{f_{y}}{f}\partial_{x}+\frac{g_{x}}{g}\partial_{y})+(\frac{f_{y}}{f})_{_{y}}\partial_{x}+\frac{g_{x}}{g}(-\frac{g_{x}g}{f^{2}}\partial_{x}+\frac{g_{y}}{g}\partial_{y})+(\frac{g_{x}}{g})_{_{y}}\partial_{y}
\\ 
&&\!\!\!\!+\left[\frac{gg_{x}}{f^{2}}(\frac{f_{x}}{f}\partial_{x}-\frac{ff_{y}}{g^{2}}\partial_{y})+(\frac{gg_{x}}{f^{2}})_{_{x}}\partial_{x}-\frac{g_{y}}{g}(\frac{f_{y}}{f}\partial_{x}+\frac{g_{x}}{g}\partial_{y})-(\frac{g_{y}}{g})_{_{x}}\partial_{y}\right]
\\
&=&\!\!\!\!\left[\frac{f_{y}^{2}}{f^{2}}-\frac{g_{x}^{2}}{f^{2}}+\frac{gf_{x}g_{x}}{f^{3}}-\frac{g_{y}f_{y}}{fg}+(\frac{f_{y}}{f})_{_{y}}+(\frac{gg_{x}}{f^{2}})_{_{x}}\right]\partial_{x}
\end{eqnarray*}
\section{Geodesic vector fields} 
\normalfont A codimension-one foliation $\mathcal{F}$ on a Riemannian manifold $M$ is \textit{totally geodesic} if their leaves are totally geodesic submanifolds of $M$. If an unit vector field 
\begin{eqnarray*}
	U_{(x,y)}=a(x,y)\partial_{x}+b(x,y)\partial_{y} 
\end{eqnarray*}
defined on $\mathbb{R}\times_{(f,g)}\mathbb{R}$ is geodesic, then its integral curves determines a totally geodesic foliation on $\mathbb{R}\times_{(f,g)}\mathbb{R}$. The vector field $U$ is geodesic $\Leftrightarrow$ $\nabla_{U}U=0$. But
\begin{eqnarray*}	\nabla_{U}U\!\!\!\!\!&=&\!\!\!a\nabla_{\partial_{x}}U+b\nabla_{\partial_{y}}U
\\	&=&\!\!\!a\left(\nabla_{\partial_{x}}(a\partial_{x}+b\partial_{y})\right)+b\left(\nabla_{\partial_{y}}(a\partial_{x}+b\partial_{y})\right)
\\
&=&\!\!\!a\left(a\nabla_{\partial_{x}}\partial_{x}+a_{x}\partial_{x}+b\nabla_{\partial_{x}}\partial_{y}+b_{x}\partial_{y}\right)\!+\!b\left(a\nabla_{\partial_{y}}\partial_{x}+a_{y}\partial_{x}+b\nabla_{\partial_{y}}\partial_{y}+b_{y}\partial_{y}\right)
\\
&=&\!\!\!a\left(a(\frac{f_{x}}{f}\partial_{x}-\frac{ff_{y}}{g^2}\partial_{y})+a_{x}\partial_{x}+b(\frac{f_{y}}{f}\partial_{x}+\frac{g_{x}}{g}\partial_{y})+b_{x}\partial_{y}\right)
\\
&&\!\!\!+b\left(a(\frac{f_{y}}{f}\partial_{x}+\frac{g_{x}}{g}\partial_{y})+a_{y}\partial_{x}+b(-\frac{gg_{x}}{f^{2}}\partial_{x}+\frac{g_{y}}{g}\partial_{y})+b_{y}\partial_{y}\right)
\\
&=&\!\!\!\left(a^{2}\frac{f_{x}}{f}-b^{2}\frac{g_{x}g}{f^{2}}+2ab\frac{f_{y}}{f}+aa_{x}+ba_{y}\right)\partial_{x}
\\
&&\!\!\!+\left(-a^{2}\frac{ff_{y}}{g^{2}}+b^{2}\frac{g_{y}}{g}+2ab\frac{g_{x}}{g}+ab_{x}+bb_{y}\right)\partial_{y}
\end{eqnarray*}
Therefore, the condition $\nabla_{U}U=0$ is equivalent to the following system
\[
\left\{
      \begin{array}{ll}
      a^{2}\frac{f_{x}}{f}-b^{2}\frac{g_{x}g}{f^{2}}+2ab\frac{f_{y}}{f}+aa_{x}+ba_{y}=0\\
      \\
      -a^{2}\frac{ff_{y}}{g^{2}}+b^{2}\frac{g_{y}}{g}+2ab\frac{g_{x}}{g}+ab_{x}+bb_{y}=0
      \end{array}
\right.
\]
\theorem Consider the following unit vector field defined on $\mathbb{R}\times_{(f,g)}\mathbb{R}$
\begin{eqnarray*}
	U_{(x,y)}:=\frac{1}{\sqrt{2}f(x,y)}\partial_{x}+\frac{1}{\sqrt{2}g(x,y)}\partial_{y}
\end{eqnarray*}
Then the vector field $U$ above is an unit geodesic vector field on $\mathbb{R}\times_{(f,g)}\mathbb{R}$ iff 
\begin{eqnarray*}
	f_{y}=g_{x}
\end{eqnarray*}
\proof In this case we have $a=1/\sqrt{2}f$ and $b=1/\sqrt{2}g$. Substituting this in the above system we obtain
\[
\left\{
      \begin{array}{ll}    \frac{1}{2f^{2}}\cdot\frac{f_{x}}{f}-\frac{1}{2g^{2}}\cdot\frac{g_{x}g}{f^{2}}+\frac{1}{fg}\cdot\frac{f_{y}}{f}+\frac{1}{\sqrt{2}f}\cdot\frac{\partial}{\partial x}(\frac{1}{\sqrt{2}f})+\frac{1}{\sqrt{2}g}\cdot\frac{\partial}{\partial y}(\frac{1}{\sqrt{2}f})=0\\
\\
-\frac{1}{2f^{2}}\cdot\frac{ff_{y}}{g^{2}}+\frac{1}{2g^{2}}\cdot\frac{g_{y}}{g}+\frac{1}{fg}\cdot\frac{g_{x}}{g}+\frac{1}{\sqrt{2}f}\cdot\frac{\partial}{\partial x}(\frac{1}{\sqrt{2}g})+\frac{1}{\sqrt{2}g}\cdot\frac{\partial}{\partial y}(\frac{1}{\sqrt{2}g})=0
      \end{array}
\right.
\]
After simplifications, both equations are equivalent to the equation $f_{y}=g_{x}$ $\square$
\\
\\
\textbf{Example.} If $f(x,y):=\mathrm{e}^{x+y}+\sin^{2}(x)+1$ and $g(x,y):=\mathrm{e}^{x+y}+1$, then $f_{y}=g_{x}$ and the following unit vector field of $\mathbb{R}\times_{(f,g)}\mathbb{R}$
\begin{eqnarray*}	U_{(x,y)}:=\frac{1}{\sqrt{2}(\mathrm{e}^{x+y}+\sin^{2}(x)+1)}\partial_{x}+\frac{1}{\sqrt{2}(\mathrm{e}^{x+y}+1)}\partial_{y}
\end{eqnarray*}
is a geodesic vector field.
\section{Foliations with geodesic normal vector field}
We recall some results about totally geodesic foliations. Let $M^{n+1}$ be an orientable Riemannian manifold and $\mathcal{F}$ be a codimension-one $\mathcal{C}^{\infty}-$foliation on $M$. Suppose that $\mathcal{F}$ is \textit{transversely orientable}, i.e., we may choose a differentiable unit vector field $N\in\mathfrak{X}(M)$ normal to the leaves of $\mathcal{F}$. 
\\
\\
Given a point $p\in M$ we may always choose an orthonormal frame field
\begin{eqnarray*}
	\left\{e_{1},e_{2},\ldots, e_{n},e_{n+1}\right\}
\end{eqnarray*}
defined in a neighborhood of $p$ and such that the vectors $e_{1},\ldots,e_{n}$ are tangent to the leaves of $\mathcal{F}$ and $e_{n+1}=N$. The \textit{Ricci curvature in the direction $e_{n+1}=N$} is
\begin{eqnarray*}
	Ric(N):=\frac{1}{n}\sum\limits_{i=1}^{n}\left\langle \mathrm{R}(e_{i},N)e_{i},N\right\rangle
\end{eqnarray*}
The \textit{divergence} of a vector field $V\in\mathfrak{X}(M)$ is locally defined as
\begin{eqnarray*}
	div(V):=\sum\limits_{k=1}^{n+1}\left\langle \nabla_{e_{k}}V,e_{k}\right\rangle
\end{eqnarray*}
In [3] the Author proved the following 

\theorem 
Let $\mathcal{F}$ be a codimension-one foliation of a complete Riemannian manifold $M^{n}$ and let $N$ be an unit vector field normal to the leaves of $\mathcal{F}$. Then
\begin{eqnarray*}
	\mathcal{F}\ \textit{is totally geodesic}\Longleftrightarrow Ric(N)-\frac{1}{n}div(\nabla_{N}N)= 0
\end{eqnarray*}
\normalfont As a consequence we obtain

\corollary Let $\mathcal{F}$ be a codimension-one foliation of a complete Riemannian manifold $M$ and let $N$ be a geodesic unit vector field normal to the leaves of $\mathcal{F}$. Then $\mathcal{F}$ is a totally geodesic foliation if, and only if, $Ric(N)=0$.
\\
\\
\normalfont The following theorem gives us some examples of such foliations

\theorem Let $f=f(x)\geq 1$ and $g=g(y)\geq 1$ be differentiable functions. Consider the following unit vector field defined on $\mathbb{R}\times_{(f,g)}\mathbb{R}$
\begin{eqnarray*}
	N_{(x,y)}:=\frac{1}{\sqrt{2}f(x)}\partial_{x}+\frac{1}{\sqrt{2}g(y)}\partial_{y}
\end{eqnarray*}
Then $N$ is a normal vector field of a totally geodesic foliation of $\mathbb{R}\times_{(f,g)}\mathbb{R}$.
\proof We have $f_{y}=g_{x}=0$ and then $N$ is a geodesic vector field. Moreover, for this choice of functions, the double-twisted product $\mathbb{R}\times_{(f,g)}\mathbb{R}$ is flat (see the expressions of curvature $\mathrm{R}$ in the $\S 2$) and complete (by Lemma 1). Therefore, $Ric(N)=0$ and (by Corollary 1) $N$ is a normal vector field of a totally geodesic foliation of $\mathbb{R}\times_{(f,g)}\mathbb{R}$ $\square$ 
\\
\\
\textbf{Example.} If $f(x):=\mathrm{e}^{x}+1$ and $g(y):=\sin^{2}(y)+1$, the following vector field
\begin{eqnarray*}	N_{(x,y)}:=\frac{1}{\sqrt{2}(\mathrm{e}^{x}+1)}\partial_{x}+\frac{1}{\sqrt{2}(\sin^{2}(y)+1)}\partial_{y}
\end{eqnarray*}
is a geodesic vector field and normal to the leaves of a totally geodesic foliation of $\mathbb{R}\times_{(f,g)}\mathbb{R}$.


\end{document}